\documentclass[oneside,reqno]{amsart}
\usepackage{amsthm}
\usepackage{amssymb}
\usepackage{amsmath,amssymb}
\usepackage{amsmath} 
\theoremstyle{plain}
\newtheorem{theorem}{Theorem}[section]

\newtheorem{lemma}[theorem]{Lemma}

\theoremstyle{definition}

\newtheorem{remark}[theorem]{Remark}

\newtheorem{problem}{Problem}

\newcommand{\C}{\mathbb{C}}

\newcommand{\cS}{\mathcal{S}}
\newcommand{\grG}{\mathcal{G}}
\newcommand{\grK}{\mathcal{K}}
\newcommand{\grA}{\mathcal{A}}
\newcommand{\grN}{\mathcal{N}}
\newcommand{\grH}{\mathcal{H}}
\newcommand{\supp}{\mathrm{supp}}

\newcommand{\cA}{\mathcal{A}}

\newcommand{\cE}{\mathcal{E}}

\newcommand{\cone}{\Omega}
\newcommand{\ip}[2]{( #1,#2)}
\newcommand{\dup}[2]{\langle #1,#2\rangle}

\title[Atomic decompositions of mixed norm Bergman spaces]{Atomic decompositions of mixed norm Bergman spaces  on tube type domains}

\author{Jens Gerlach Christensen} 
\address{ Department of Mathematics, Colgate University, 13 Oak Drive, Hamilton NY 13346} 
\email{jchristensen@colgate.edu}
\urladdr{http://www.math.colgate.edu/~jchristensen}

\begin{document}
\begin{abstract}
  We use the author's previous work on atomic decompositions 
  of Besov spaces with spectrum on symmetric cones, to derive
  new atomic decompositions for Bergman spaces on 
  tube type domains. It is related to work by
  Ricci and Taibleson who derived decompositions
  for classical Besov spaces from atomic decompositions of
  Bergman spaces on the upper half plane.
  Moreover, for this class of domains our method is an alternative to 
  classical results by Coifman and Rochberg,
  and it works for a larger range of Bergman spaces.
\end{abstract}

\maketitle

\section{Introduction}
In this paper we suggest a new approach to 
atomic decompositions for the Bergman spaces
on tube type domains. Atomic decompositions have
previously been obtained for the upper half plane (unit disc)
in \cite{Ricci1983,Christensen2009,Pap2012},  and for the unit ball 
\cite{Zhu2005,Christensen2017} as well as
other bounded symmetric domains 
\cite{Coifman1980,Luecking1985}.
These decompositions have typically been attained by investigating
oscillations of the Bergman kernel and in most cases the atoms are samples
of the Bergman kernel. In the case of the unit ball
such oscillations can be estimated for all parameters,
but for higher rank spaces it is more complicated. The issue is connected
to the question of boundedness of the 
Bergman projection on bounded symmetric domains, which is still an
open problem. 
Advances to answer this problem have recently been made
in the case of tube type domains in \cite{Bekolle2004,Debertol2005} and
for general domains in \cite{Nana2013}. 
In the case of tube type domains over forward light cones the 
problem is now solved \cite{Bonami2015}.

In this paper we will concentrate on tube type domains, and we
use Fourier-Laplace extensions to transfer known atomic 
decompositions for Besov spaces \cite{Christensen2013a} to Bergman spaces.
This allows us
to narrow the gap in the atomic decompositions from \cite{Coifman1980}
in the case of tube type domains.
Moreover, these decompositions are for mixed norm Bergman spaces
which have not previously been dealt with. 
We would like to mention that this work seems to be in the reverse direction
of the paper \cite{Ricci1983} which uses atomic decompositions for
the mixed norm Bergman spaces on the upper half plane to get
atomic decompositions for the Besov spaces. See also \cite{Feldman1991}
for related work for Bergman spaces on the unit ball in relation to
Besov spaces on the Heisenberg group.

\section{Symmetric cones}
For an introduction to symmetric cones we refer to the book
\cite{Faraut1994}.
Let $V$ be a Euclidean vector space over the real numbers of finite
dimension $n$.
A subset $\cone$ of $V$ is a cone if $\lambda\cone \subseteq \cone$
for all $\lambda>0$. Assume $\cone$ is open and convex, and
define the open dual cone $\cone^*$ by
\begin{equation*}
  \cone^* = \{y\in V \mid \ip{x}{y} > 0\,\, 
  \text{for all non-zero $x\in \overline{\cone}$} \}.
\end{equation*}
The cone $\cone$ is called symmetric if $\cone=\cone^*$, and
the automorphism group
\begin{equation*}
  \grG(\cone) = \{ g\in\mathrm{GL}(V)\mid g\cone=\cone   \}
\end{equation*}
acts transitively on $\cone$. 
Notice that the group $\grG(\cone)$ is semisimple.
Define the characteristic function of $\cone$ by
\begin{equation*}
  \varphi(x) = \int_{\cone^*} e^{-\ip{x}{y}}\,dy,
\end{equation*}
then 
\begin{equation*}
  \varphi(gx) = |\det(g)|^{-1}\varphi(x).
\end{equation*}
Also,
\begin{equation}
  \label{eq:invintcone}
  f\mapsto \int_\cone f(x)\varphi(x)\,dx
\end{equation}
defines a $\grG(\cone)$-invariant measure on $\cone$.
The connected component $\grG_0(\cone)$ of $\grG(\cone)$
has Iwasawa decomposition
\begin{equation*}
  \grG_0(\cone)= \grK \grA \grN
\end{equation*}
where $\grK = \grG_0(\cone)\cap \mathrm{O}(V)$ is compact, 
$\grA$ is abelian and $\grN$ is nilpotent.
The unique fixed
point in $\cone$ for the mapping $x\mapsto \nabla \log\varphi(x)$
will be denoted $e$, and we note that $\grK$ fixes $e$.
The connected solvable subgroup $\grH = \grA \grN$ of $\grG_0(\cone)$ 
acts simply transitively on $\cone$ and the integral
(\ref{eq:invintcone}) thus also
defines the left-Haar measure on $\grH$.
Denote by $\mathcal{S}(V)$ the space of rapidly decreasing
smooth functions with topology induced by the semi-norms
\begin{equation*}
  \| f\|_{k} 
  = \sup_{|\alpha|\leq k} \sup_{x\in V} |\partial^{\alpha} f (x)|(1+|x|)^k.
\end{equation*}
Here $\alpha$ is a multi-index, $\partial^\alpha$ denotes usual partial
derivatives of functions,  and $k\geq 0$ is an integer.
For $f\in \mathcal{S}(V)$ the Fourier transform is defined by
\begin{equation*}
  \widehat{f}(w) = \frac{1}{(2\pi)^{n/2}}\int_V f(x)e^{-i\ip{x}{w}}\,dx
  \text{\ for $w\in V$}.
\end{equation*}
The convolution 
$$
f*g(x) = \int_V f(y)g(x-y)\,dy
$$
of functions $f,g\in \mathcal{S}(V)$ satisfies
\begin{equation*}
  \widehat{f*g}(w) = \widehat{f}(w)\widehat{g}(w).
\end{equation*}
The space $\mathcal{S}'(V)$ of tempered distributions
is the linear dual of $\mathcal{S}(V)$.
For functions on $V$ define
$\tau_x f(y) = f(y-x)$, $f^\vee(y) = f(-y)$ and
$f^*(y)=\overline{f(-y)}$.
Convolution of $f\in \mathcal{S}'(V)$ and $\phi\in\mathcal{S}(V)$
is defined by 
$$
f*\phi(x) = f(\tau_x\phi^\vee).
$$
As usual, the Fourier transform extends to tempered distributions
by duality.
The space of rapidly decreasing smooth
functions with Fourier transform vanishing on $\cone$ is
denoted $\mathcal{S}_\cone$. It is a closed subspace of $\mathcal{S}(V)$
and will be equipped with the subspace topology.

The space $V$ can be equipped with a Jordan algebra structure
such that $\overline{\cone}$ is identified with the set of all
squares. This gives rise to the notion of a determinant
$\Delta(x)$. We only need the fact that the determinant
is related to the characteristic function $\varphi$ by
\begin{equation*}
  \varphi(x) = \varphi(e)\Delta(x)^{-n/R},
\end{equation*}
where $R$ denotes the rank of the cone.
If $x=ge$ we have
\begin{equation} \label{eq:determinantrelation}
  \Delta(x) = \Delta(ge) = |\mathrm{Det}(g)|^{R/n}.
\end{equation}

\section{Besov spaces related to symmetric cones}
The cone $\cone$ can be identified as a Riemannian manifold
$\cone = \grG_0(\cone)/\grK$ where $\grK$ is the compact
group fixing $e$. The Riemannian metric
in this case is defined by
\begin{equation*}
  \dup{u}{v}_y = \ip{g^{-1}u}{g^{-1}v}
\end{equation*}
for $u,v$ tangent vectors to $\cone$ at $y=ge$. Denote the balls
of radius $\delta$ centered at $x$ by $B_\delta(x)$. For $\delta>0$ and 
$\lambda\geq 2$ the points $\{ x_j\}$ are called a 
$(\delta,\lambda)$-lattice if
\begin{enumerate}
\item $\{ B_\delta(x_j)\}$ are disjoint, and
\item $\{ B_{\lambda\delta}(x_j)\}$ cover $\cone$.
\end{enumerate}

We now fix a $(\delta,\lambda)$-lattice $\{ x_j\}$ with $\delta=1/2$ and 
$\lambda=2$.
Then there are functions $\psi_j\in\mathcal{S}_\cone$,
such that $0\leq\widehat{\psi}_j\leq 1$, $\mathrm{supp}(\widehat{\psi}_j) \subseteq B_2(x_j)$,
$\widehat{\psi}_j$ is one on $B_{1/2}(x_j)$ and $\sum_j \widehat{\psi}_j =1$ on $\cone$.
Using this decomposition of the cone, the Besov space norm for $1\leq p,q<\infty$ and $\nu\in \mathbb{R}$ is
defined in \cite{Bekolle2004} by
\begin{equation*}
  \| f\|_{B^{p,q}_\nu} = \left( \sum_j \Delta(x_j)^{-\nu}\| f*\psi_j\|_p^q \right)^{1/q}.
\end{equation*}
The Besov space $B^{p,q}_\nu$ consists of 
the equivalence classes of tempered distributions $f$ in
$(\mathcal{S}_\cone)' \simeq 
\{ f\in \mathcal{S}'(V) \mid \supp(\widehat{f})\subseteq \overline{\cone}\}
/\mathcal{S}'_{\partial \cone}$ 
for which $\| f\|_{B^{p,q}_s} <\infty$.

Define the index
$$
\widetilde{q}_{\nu,p} = \frac{\nu+n/R-1}{(n/Rp')-1}
$$
if $n/R>p'$ and
set $\widetilde{q}_{\nu,p} = \infty$ if $n/R\leq p'$.
The following results from \cite{Bekolle2004} states when the Besov spaces
are included in the space of tempered distributions $\cS'(V)$.
\begin{lemma}
  Let $\nu>0$, $1\leq p<\infty$ and $1\leq q < \widetilde{q}_{\nu,p}$.
  Then for every $f\in B^{p,q}_\nu$ the series $\sum_j f*\psi_j$ converges
  in the space $\cS'(V)$, and the correspondence
  $$
  f+\cS'_{\partial \cone} \mapsto f^\sharp = \sum_j f*\psi_j
  $$
  is continuous, injective and does not depend on the particular
  choice of $\{ \psi_j\}$.
\end{lemma}

The main result from \cite{Christensen2013a} is that
the quasiregular representation of the group $\grH \rtimes V$ 
can be used to obtain atomic decompositions for these Besov spaces.
We summarize the result in
\begin{theorem}
  Let $\psi$ in $\cS_\cone$ be such that $\widehat{\psi}$ is 
  compactly supported. Then there is an index set $I$, 
  a set $\{ (h_i,x_i)\}_{i\in I}\subseteq \grH\times V$, a Banach sequence space
  $b^{p,q}_\nu(I)$, a set of continuous functionals
  $\{ c_i:B^{p,q}_\nu\to\mathbb{C}\}_{i\in I}$, and a constant $C>0$ such that
  
  \begin{enumerate}
  \item $f(x) = \sum_{i\in I} c_i(f) 
    \frac{1}{\sqrt{\det(h_i)}}\psi (h_i^{-1}(x-x_i))$
    with convergence in norm in $B^{p,q}_\nu$
  \item $\| c_i(f)\|_{b^{p,q}_\nu} \leq C \| f\|_{B^{p,q}_\nu}$
  \item if     $\{ \lambda_i\} \in b^{p,q}_\nu$ then
    $$f(x)= \sum_i \lambda_i  \frac{1}{\sqrt{\det(h_i)}}\psi (h_i^{-1}(x-x_i))$$
    is in $B^{p,q}_\nu$ and 
    $\| f\|_{B^{p,q}_\nu} \leq C \| \{ \lambda_i\}_{b^{p,q}_\nu}\|$.
  \end{enumerate}
\end{theorem}

\begin{remark}
  At this stage it is appropriate to describe the sequence of points
  $\{ (h_i,x_i)\}_{i\in I}$ and the sequence
  space $b^{p,q}_\nu$ in some detail. 
  For this one chooses a covering $\{ U_i\}_{i\in I}$
  of the space $\cone\times V$. This covering is chosen such
  that each $U_i$ is a translate of a fixed relatively compact neighbourhood
  $U$ of $\{ e\}\times \{ 0\}$
  by some element $(h_i,x_i)$ of the semidirect product
  $\grH\rtimes V$. Moreover, the sets $U_i$ have the finite 
  overlapping property, that is, there is an $N$ such that each 
  set $U_i$ overlap at most $N$ others.
  A sequence $\{ \lambda_i\}$ is in $b^{p,q}_\nu$ if
  $$
   \| \{\lambda_i\}\|_{b^{p,q}_\nu} 
   :=\left(\int_\cone \left( \int_V \sum_{i\in I} 
     |\lambda_i| 1_{U_i}(x,t)^{p}\,dt\right)^{q/p}
   \Delta(x)^{\nu - qn/(2R)-n/R} \,dx\right)^{1/q}
  $$
  is finite. If $p=q$ this is a $\Delta^{\nu - pn/(2R)-n/R}$-weighted 
  $\ell^p$-space.
\end{remark}

\section{Bergman spaces on tube type domains.}

In this section we introduce the Bergman spaces on the
tube type domains, and describe the isomorphism between
a range of Besov spaces and Bergman spaces.

Let $T=\{ z = x+iy \mid x\in V,y\in\cone \}$ be the tube type
domain related to the symmetric cone $\cone$. 
For $1\leq p,q<\infty$ and $\nu>0$
define the weighted Lebesgue space $L^{p,q}_\nu$ on the tube
type domain to consist of the equivalence classes of
measurable functions on $\cone$ for which the norm
$$
\| F\|_{L^{p,q}_\nu} := \left(\int_\cone \left(\int_V |F(x+iy)|^p\,dx \right)^{q/p} 
\Delta^{\nu-n/r}(y)\,dy\right)^{1/q}
$$
is finite.
Here $dx$ and $dy$ denote the usual Lebesgue measures on $V$ and $\cone$.
The mixed norm Bergman space $\cA^{p,q}_\nu$ on $\Omega$
consists of the holomorphic functions in $L^{p,q}_\nu$.
It is well-known that this is a reproducing kernel Banach space,
that is, for every $z\in T$  the mapping
$F\mapsto F(z)$ is continuous from $\cA^{p,q}_\nu$ to $\C$.

The special case $p=q=2$ and $\nu=n/R$ is the usual Bergman space
and the reproducing kernel in this case is
$$
B(z,w) = B_{n/r}(z,w) = c(\nu) \Delta\left( \frac{z-\overline{w}}{i} \right)^{-2n/r},
$$
which will be called the Bergman kernel.

Following \cite{Bekolle2004} we now define the Fourier-Laplace extensions
of elements in the Besov spaces. This extension only works for
Besov spaces which can be naturally imbedded 
in the usual space of tempered distributions.
This introduces a restriction in the range of indices that can be used.
We summarize the results from \cite{Bekolle2004} that we need.

Define the Fourier-Laplace extension of a tempered distribution $f$
whose Fourier transform $\widehat{f}$ is supported on $\overline{\cone}$ by
$$
\cE f = \int \widehat{f}(w)e^{iz\cdot w} \,dw
$$
for $z\in T$.
For $1\leq q < \widetilde{q}_{\nu,p}$ the Besov space can be indentified
with a space of such distributions, and therefore
we can define
$$
\widetilde{\cE} f = \cE \sum_j f*\psi_j
$$
for $f\in B^{p,q}_\nu$.
Define the index
$$
q_{\nu,p} = \min(p,p') \frac{\nu+n/R-1}{n/R-1}
$$
when $n>R$ and set $q_{\nu,p}=\infty$ when $n=1$.
Notice that $2<q_{\nu,p}\leq \widetilde{q}_{\nu,p}$ so when 
$1\leq q < q_{\nu,p}$ all elements
of the Besov space $B^{p,q}_\nu$  can be identified with
tempered distributions whose Fourier-Laplace extensions
are in the Bergman space $\cA^{p,q}_\nu$.
We have
\begin{theorem}
  If $\nu> n/R-1$, $1\leq p<\infty$ 
  and $1\leq q<q_{\nu,p}$, then 
  the mapping
  $\widetilde{\cE}: B^{p,q}_\nu \to \cA^{p,q}_\nu$ is
  an isomorphism.
  Moreover,
  $$
  \lim_{y\to 0} F(x+iy) =f(x)
  $$
  in both $\cS'(V)$ and $B^{p,q}_\nu$.
\end{theorem}

\section{Atomic decomposition of Bergman spaces}

In this section we merge the results from \cite{Christensen2013a}
and \cite{Bekolle2004} to obtain atomic decompositions for
Bergman spaces on the tube type domains. This will give an alternate 
approach to the atomic decompositions found in \cite{Coifman1980}.
This new approach allows for a large class of atoms, and moreover,
through the Paley-Wiener theorem the decay properties of these
atoms are quite well known. Note, that the decomposition
from \cite{Coifman1980} uses samples of the Bergman kernel,
but the Bergman kernel is not among the possible atoms
with the new approach. The paper \cite{Ricci1983} goes in the
opposite direction and uses \cite{Coifman1980} to obtain
atomic decompositions the Besov spaces. It would of course be interesting
to investigate how to completely align the two methods
and to determine exactly which atoms can be moved from the
Besov setting and to the Bergman setting. In this paper we
are clearly only dealing with a subset of possible atoms.

Let $F\in \cA^{p,q}_\nu$ then $f=\widetilde{\cE}^{-1}F$ is in
$B^{p,q}_\nu$ with equivalent norms 
and can be decomposed as
$$
f(x) = \sum_{i\in I} c_i(f) 
    \frac{1}{\sqrt{\det(h_i)}}\psi (h_i^{-1}(x-x_i)).
$$
Since $\psi\in \cS_0$ with compactly supported
Fourier transform is in every Besov space, we get
$$
\psi_i (z) := \widetilde{\cE}(\psi (h_i^{-1}(\cdot-x_i)))(z) 
= \widetilde{\cE}\psi(h^{-1}(z-x_i)).
$$
This results in the following atomic decompositions.
\begin{theorem}
  \label{thm:maintheorem}
  Let $\nu > n/R-1$, $1\leq p <\infty$ and $1\leq q < q_{\nu,p}$.
  There is a sequence $\{ d_i\}$ of functionals and atoms $\{ \psi_i \}$
  parameterized by appropriate $\{ (h_i,x_i)\} \subseteq \grH\rtimes V$ 
  such that
  \begin{enumerate}
  \item if $F\in \cA^{p,q}_\nu$ then
    $$
    F(z) = \sum_{i\in I} d_i(F) 
    \frac{1}{\sqrt{\det(h_i)}} \psi_i(z).
    $$
    where $d_i(F)=c_i(\widetilde{\cE}^{-1}F)$ satisfy
    that $\| \{ d_i(F)\}\|_{b^{p,q}_\nu} \leq C \| F\|_{\cA^{p,q}_\nu}$,
    \item and
      if $\{d_i\} \in b^{p,q}_\nu$, then
      $$
      F(z) = \sum_{i\in I} d_i 
      \frac{1}{\sqrt{\det(h_i)}} \psi_i(z).
      $$
      is in $\cA^{p,q}_\nu$ and 
      $\| F\|_{\cA^{p,q}_\nu} \leq C \| \{ d_i\}\|_{b^{p,q}_\nu}$.
    \end{enumerate}
\end{theorem}

\begin{remark}
Notice that for the case of cones of rank 2 (for example the forward light cones)
this theorem can be extended to the larger range $1\leq q < \widetilde{q}_{\nu,p}$
which is the entire range of $q$ for which Laplace extensions can be defined. See \cite{Bonami2015}.

\end{remark}

\section{Comparison with previous results 
  and some open problems}

To demonstrate how this work extendeds the range of Bergman spaces
for which atomic decompositions can be found, we 
now compare our atomic decompositions to classical results
due to Coifman and Rochberg 
\cite{Coifman1980} in the special case of tube type domains.
Their results only work for rank one spaces, but with minor
modifications this issue can be addressed via
the Forelli-Rudin estimates from
Theorem 4.1 in \cite{Faraut1990} or 
Corollary II.4 in \cite{Bekolle1995a}.
The latter result was used in \cite{Bekolle1998}
to correct and generalize the 
atomic decompositions of Coifman and Rochberg to also include
two non-symmetric domains.

We first summarize the atomic decompsitions from \cite{Coifman1980,Bekolle1998}.
Let $V$ be an open convex cone in $\mathbb{R}^m$ and let 
$F$ be a $V$-valued Hermitian form on $\mathbb{C}^n$.
The open subset of 
of elements $(z,w)$ in $\mathbb{C}^m\times \mathbb{C}^n$ 
for which $\mathrm{Im}(z)-F(w,w)\in V$ is called a
Siegel domain of type II. The domain is called symmetric
if it is also a symmetric space.
Let $D$ be a symmetric Siegel domain of type II
and let $B$ denote the associated Bergman kernel.
Coifman and Rochberg use a parametrization of Bergman spaces
that differs from the one used earlier in this paper.
In their notation the Bergman space $\widetilde{\cA}^p_r$
consists of holomorphic functions
for which the following norm is bounded
$$
\| F\|_{\widetilde{\cA}^{p}_r}
= \int_D |F(z)|^p B(z,z)^{-r} dz.
$$

\begin{theorem}
  \label{thm:coifmanrochberg}
  Let $p\geq 1$ and assume that
  $-\epsilon_D+\gamma_D(p-1) < r < \infty$.
  Given $\theta > p(1-\epsilon_D)+\epsilon_D+\gamma_D-2-r$ there is
  a lattice $\{ \xi_i\}$ in $D$ and a constant $C>0$ 
  such that for $F$ in $\widetilde{\cA}^p_r$  we have 
  $$
  F(z) 
  = \sum_i \lambda_i(F) 
  \left(\frac{B(z,\xi_i)^2}{B(\xi_i,\xi_i)}\right)^{\frac{1+r}{p}}
  \left( \frac{B(z,\xi_i)}{B(\xi_i,\xi_i)} \right)^{\frac{\theta}{p}},
  $$
  and $\sum_i |\lambda_i(F)|^p \leq C \| F\|_{\widetilde{\cA}^p_r}^p$.
  Moreover, if $\{ \lambda_i\}\in \ell^p$ then the series
$$
  F(z) 
  = \sum_i \lambda_i
  \left(\frac{B(z,\xi_i)^2}{B(\xi_i,\xi_i)}\right)^{\frac{1+r}{p}}
  \left( \frac{B(z,\xi_i)}{B(\xi_i,\xi_i)} \right)^{\frac{\theta}{p}}.
  $$
  defines a function in $\widetilde{\cA}^p_r$ and $\| F\|_{\cA^p_r} \leq C \| \{ \lambda_i\}\|_{\ell^p}.$
\end{theorem}

\begin{remark}
  The constants are given by
  $\epsilon_D=1/G$ and $\gamma_D=(R-1)a/(2G)$, where $G$ is the genus
  and $a$ is another structural constant. See, for example, \cite{Faraut1990}
  for a full explanation of these constants.
\end{remark}

We will now establish the range of $r$ which work for tube type domains.
The connection between the number 
$\nu$ from Theorem \ref{thm:maintheorem} in the special case of $p=q$
and the number $r$ from Theorem \ref{thm:coifmanrochberg}
is 
$$
\nu = \frac{2nr}{R}+ \frac{n}{R}.
$$
Moreover, the structural constants $\epsilon_D$ and $\gamma_D$
for tube type domains are
$$
\epsilon_D = \frac{R}{2n} \qquad \text{and} \qquad \gamma_D = \frac{1}{2} - \epsilon_D.
$$
With this in mind the conditions
$\nu > \frac{n}{R}$ and $p<\frac{\nu+n/R-1}{n/R-1}+1$ from Theorem~\ref{thm:maintheorem} rewrite
into
$$
r> \max\{ -\epsilon_D, -\frac{3}{2}+p(1-\epsilon_D)+\epsilon_D-\frac{p}{2} \}.
$$
The result from Theorem \ref{thm:coifmanrochberg} works for
$$
r> \max\{ -\epsilon_D+\gamma_D, -\frac{3}{2}+p(1-\epsilon_D) \},
$$
when restricted to the case $\theta=0$.
Since for tube type domains over cones $\epsilon_D = R/(2n) \leq 1/2$
and $p\geq 1$ we see that the atomic decompositions in
Theorem ~\ref{thm:maintheorem} work for a larger range than those of Theorem ~\ref{thm:coifmanrochberg}.

\begin{remark}
  The reason we restrict to the 
  $\theta=0$ when comparing the two methods is, that 
  the the atomic decompositions provided by Laplace extensions
  is connected to the discrete series representation of
  the automorphism group on the tube type domain, and 
  we might as well transfer this result to 
  the bounded realization of the domain. Therefore Theorem \ref{thm:maintheorem}
  can be transfered to Bergman spaces on the bounded symmetric domain.
  It thus makes sense to compare with the version of
  Theorem \ref{thm:coifmanrochberg} which also can be transfered to
  the bounded realization, that is, $\theta=0$.
\end{remark}

We finish this section with a list of open problems connected
to the results of this paper.

\begin{problem}\label{problem1}
  The atoms from Theorem~\ref{thm:maintheorem} 
  do not include samples of the Bergman kernel
  as in \cite{Coifman1980,Bekolle1998}. The reason is that the atoms we use
  are extensions of compactly supported smooth functions, and 
  by the Paley-Wiener Theorem these cannot include the Bergman kernel.
  This means that Theorem~\ref{thm:maintheorem} does not
  include Theorem~\ref{thm:coifmanrochberg} as a special case.
  We believe that it is possible to overcome this issue.
  The atomic decompositions
  in \cite{Christensen2013a} build on irreducible, unitary, and 
  integrable group 
  representations, and therefore a much larger class of atoms for the 
  Besov spaces can
  be used via \cite{Feichtinger1989}. It would be interesting to
  see if the Laplace extensions of this larger class of atoms for the Besov
  spaces would include the Bergman kernel in order to
  obtain Theorem \ref{thm:coifmanrochberg} as a special case of Theorem \ref{thm:maintheorem}.
  This question would be of interest even on the upper half plane.

\end{problem}

\begin{remark}
    Upon completion of this work the author was made aware of a related
  paper by D. B\'ekoll\'e, J. Gonessa and C. Nana \cite{Bekolle2017}. They obtain
  atomic decompositions for the exact same range of Bergman spaces. 
  In their work the atoms are 
  indeed samples of the Bergman kernel. Problem ~\ref{problem1}
  thus formulates one possible approach to uncovering the connection
  between their result and the present paper.
\end{remark}

\begin{problem}
  Another way to get a larger set of atoms including samples of 
  the Bergman kernel
  is to apply the coorbit
  theory \cite{Christensen2011} as has been done for the unit ball 
  in \cite{Christensen2017,Christensen2017a}. 
  The integral operator with positive kernel 
  from Theorem II.7 in \cite{Bekolle1995a}
  which was used to derive Theorem~\ref{thm:coifmanrochberg} 
  would play a crucial role in this approach
  (as it did on the unit ball \cite{Christensen2017,Christensen2017a}), 
  so we predict this approach would work for the same range of parameters as 
  Theorem~\ref{thm:coifmanrochberg}. It would be interesting to see
  if the use of Theorem II.7 in \cite{Bekolle1995a}
  could be avoided or refined in the context of coorbits in order to get to the same range of Bergman
  spaces as in Theorem~\ref{thm:maintheorem}.
\end{problem}

\begin{problem}
  The approach highlighted in this paper could be used in the setting
  of the unit ball by using Laplace extensions of Besov
  spaces on the Heisenberg group mentioned in \cite{Feldman1991}.
  It also seems possible to extend our approach to all bounded symmetric domains
  via work in \cite{Nana2013}.
\end{problem}

\end{document}